\documentclass[12pt]{article}
\usepackage{amsmath}
\usepackage{amsthm}
\usepackage{amsfonts}
\usepackage{amssymb}
\usepackage{graphicx}
\usepackage{xypic}
\usepackage{a4wide}
\usepackage{enumerate}
\usepackage{mathrsfs}
\newtheorem{theorem}{Theorem}[section]
\newtheorem{proposition}[theorem]{Proposition}
\newtheorem{corollary}[theorem]{Corollary}
\newtheorem{lemma}[theorem]{Lemma}
\theoremstyle{definition}

\newtheorem{definition}[theorem]{Definition}
\newtheorem{example}[theorem]{Example}

\newtheorem{remark}[theorem]{Remark}

\newcommand{\bcomod}[2]{{}^{#1}\mathcal{M}^{#2}}
\newcommand{\cat}[1]{\mathbf{#1}}
\newcommand{\cohom}[3]{\mathrm{h}_{#1}(#2,#3)}
\newcommand{\coring}[1]{\mathfrak{#1}}
\newcommand{\cotensor}[1]{\square_{#1}}
\renewcommand{\hom}[3]{\mathrm{Hom}_{#1}(#2,#3)}
\newcommand{\lcomod}[1]{{}^{#1}\mathcal{M}}
\newcommand{\ldual}[1]{{ }^{\ast}#1}
\newcommand{\rdual}[1]{#1^{\ast}}
\newcommand{\lDual}{{ }^{\ast}(-)}
\newcommand{\rDual}{ (-)^{\ast}}
\newcommand{\rcomod}[1]{\mathcal{M}^{#1}}
\newcommand{\rmod}[1]{\mathcal{M}_{#1}}
\newcommand{\lmod}[1]{{}_{#1}\mathcal{M}}
\newcommand{\tensor}[1]{\otimes_{#1}}
\newcommand{\nat}[2]{\mathrm{Nat}(#1,#2)}

\date{\empty}

\begin{document}

\title{Frobenius Functors for Corings}

\author{J. G\'omez-Torrecillas\footnote{Partially supported by the grant BFM2001-3141 from the Ministerio de Ciencia y Tecnolog\'{\i}a} \\
\normalsize Departamento de \'{A}lgebra \\
\normalsize Facultad de Ciencias \\
\normalsize Universidad de Granada\\
\normalsize E18071 Granada, Spain \\
\normalsize e-mail: \textsf{gomezj@ugr.es} \and M.~Zarouali Darkaoui\\
\normalsize Departamento de \'{A}lgebra \\
\normalsize Facultad de Ciencias \\
\normalsize Universidad de Granada\\
\normalsize E18071 Granada, Spain } \maketitle

\section*{Introduction}

Corings were introduced by M.~Sweedler in \cite{Sweedler:1975} as
a generalization of coalgebras over commutative rings to the case
of non-commutative rings, to give a formulation of a predual of
the Jacobson-Bourbaki's theorem for intermediate extensions of
division ring extensions. Thus, a coring over an associative ring
with unit is a coalgebra in the monoidal category of all
$A$--bimodules. Recently, motivated by an observation of M.
Takeuchi, namely that an entwining structure (resp. an entwined
module) can be viewed as a suitable coring (resp. as a comodule
over a suitable coring), T.~Brzezi\'{n}ski, has given in
\cite{Brzezinski:2002} some new examples and general properties of
corings. Among them, a study of Frobenius corings is developed,
extending previous results on entwining structures
\cite{Brzezinski:2000} and relative Hopf modules
\cite{Caenepeel/Militaru/Zhu:1997}.

A motivation to the study of Frobenius functors is the observation
of K.~Morita \cite{Morita:1965} that a ring extension is Frobenius
if and only the restriction of scalars functor has isomorphic left
and right adjoint functors. Thus a pair of functors $(F,G)$ is
said to be a \emph{Frobenius pair}
\cite{Castano/Gomez/Nastasescu:1999} if $G$ is left and right
adjoint to $F$ at the same time. The functors $F$ and $G$ then
known as \emph{Frobenius functors}
\cite{Caenepeel/Militaru/Zhu:1997}. These names are now standard,
and we still using them rather than the original Morita's
denomination as \emph{strongly adjoint pairs}.

\emph{Frobenius corings} (i.e., corings for which the functor
forgetting the coaction is Frobenius), have been intensively
studied in
\cite{Brzezinski:2002,Brzezinski:2003,Brzezinski/Gomez:2003,Caenepeel/DeGroot/Militaru:2002}.
In this paper we investigate Frobenius pairs between categories of
comodules over rather general corings. Precedents for coalgebras
over fields are contained in \cite{Castano/Gomez/Nastasescu:1999}
and \cite{Caenepeel/DeGroot/Militaru:2002}. We particularize to
the case of the adjoint pair of functors associated to a morphism
of corings over different base rings \cite{Gomez:2002}, which
leads to a reasonable notion of Frobenius coring extension. When
applied to corings stemming from entwining structures, we obtain
new results in this setting.

\section{Basic notations}

Throughout this paper, $A,$ $A^{\prime},$ $A^{\prime\prime},$ and
$B$ denote associative and unitary algebras over a commutative
ring $k$, and except mention clarifies opposite, $\coring{C},$
$\coring{C}^{\prime},$ $\coring{C}^{\prime\prime},$ and
$\coring{D}$ denote corings over $A,$ $A^{\prime},$
$A^{\prime\prime},$ and $B$, respectively. We recall from
\cite{Sweedler:1975} that an $A$--\emph{coring} consists of an
$A$-bimodule $\coring{C}$ with two $A$--bimodule maps
\[
\Delta : \coring{C} \longrightarrow \coring{C} \tensor{A}
\coring{C}, \qquad \epsilon : \coring{C} \longrightarrow A
\]
such that $(\coring{C} \tensor{A} \Delta) \circ \Delta = (\Delta
\tensor{A} \coring{C}) \circ \Delta$ and $(\epsilon \tensor{A}
\coring{C}) \circ \Delta = ( \coring{C} \tensor{A} \epsilon) \circ
\Delta = 1_\coring{C}$. A \emph{right}
$\coring{C}$--\emph{comodule} is a pair $(M,\rho_M)$ consisting of
a right $A$--module $M$ and an $A$--linear map $\rho_M: M
\rightarrow M \tensor{A} \coring{C}$ (the coaction) satisfying $(M
\tensor{A} \Delta) \circ \rho_M = (\rho_M \tensor{A} \coring{C})
\circ \rho_M$, $(M \tensor{A} \epsilon) \circ \rho_M = 1_M$. A
\emph{morphism} of right $\coring{C}$--comodules $(M,\rho_M)$ and
$(N,\rho_N)$ is a right $A$--linear map $f: M \rightarrow N$ such
that $(f \tensor{A} \coring{C}) \circ \rho_M = \rho_N \circ f$;
the $k$--module of all such morphisms will be denoted by
$\hom{\coring{C}}{M}{N}$. The right $\coring{C}$--comodules
together with their morphisms form the additive category
$\rcomod{\coring{C}}$. Coproducts and cokernels in
$\rcomod{\coring{C}}$ do exist and can be already computed in
$\rmod{A}$, the category of right $A$--modules. Therefore,
$\rcomod{\coring{C}}$ has arbitrary inductive limits. If
${}_A\coring{C}$ is flat, then $\rcomod{\coring{C}}$ is a
Grothendieck category. The converse is not true, as \cite[Example
1.1]{ElKaoutit/Gomez/Lobillo:2001unp} shows. When $\coring{C} = A$
with the trivial $A$--coring structure, then $\rcomod{A} =
\rmod{A}$.

Now assume that the $A'-A$--bimodule $M$ is also a left comodule
over an $A'$--coring $\coring{C}'$ with structure map $\lambda_M :
M \rightarrow \coring{C}' \tensor{A} M$. It is clear that $\rho_M
: M \rightarrow M \tensor{A} \coring{C}$ is a morphism of left
$\coring{C}'$--comodules if and only if $\lambda_M : M \rightarrow
\coring{C}' \tensor{A'} M$ is a morphism of right
$\coring{C}$--comodules. In this case, we say that $M$ is a
$\coring{C}'-\coring{C}$--bicomodule. The
$\coring{C}'-\coring{C}$--bicomodules are the objects of a
category $\bcomod{\coring{C}'}{\coring{C}}$, whose morphisms are
defined in the obvious way.

Let $Z$ be a left $A$--module and $f: X \rightarrow Y$ a morphism
in $\rmod{A}$. Following \cite[40.13]{Brzezinski/Wisbauer:2003} we
say that $f$ is \emph{$Z$--pure} when the functor $-\tensor{A} Z$
preserves the kernel of $f$. If $f$ is $Z$--pure for every $Z \in
\lmod{A}$ then we say simply that $f$ is \emph{pure} in
$\rmod{A}$.

\section{Frobenius functors between categories of
comodules}\label{Frobeniusgeneral}

Let $T$ be a $k$-algebra, and $M\in{}^{T}\mathcal{M}^{A}$. Let
$\varphi:T\longrightarrow \operatorname{End}_{A}(M)$ the morphism
of $k$--algebras given by the right $T$-module structure of the
bimodule $_{T}M_{A}$. Now, suppose moreover that
$M\in\mathcal{M}^{\coring{C}}$. Then
$\operatorname{End}_{\coring{C}}(M)$ is a subalgebra of
$\operatorname{End}_{A}(M).$ We have that $\varphi(T)\subset
\operatorname{End}_{\coring{C}}(M)$ if and only if $\rho_{M}$ is
$T$-linear.
Hence the left $T$-module structure of a $T-\coring{C}$-bicomodule
$M$ can be described as a morphism of $k$-algebras$\
\varphi:T\rightarrow \operatorname{End}_{\coring{C}}(M) $. Given a
$k$--linear functor
$F:\mathcal{M}^{\coring{C}}\rightarrow\mathcal{M}
^{\coring{D}}$, and $M\in{}^{T}\mathcal{M}%
^{\coring{C}}$, the algebra morphism
$T\overset{\varphi}{\rightarrow
}\operatorname{End}_{\coring{C}}(M)\overset{F(-)}%
{\rightarrow}End_{\coring{D}}(F(M)) $ defines a
$T-\coring{D}$-bicomodule structure on $F(M)$. We have then two
$k$-linear functors
\[
-\otimes_{T}F(-),\,F(-\otimes_{T}-):\mathcal{M}%
^{T}\times{}^{T}\mathcal{M}^{\coring{C}}\rightarrow\mathcal{M}%
^{\coring{D}}.
\]
Let $\Upsilon_{T,M}$ be the unique isomorphism of
$\coring{D}$-comodules making the following diagram commutative
\[%
\xymatrix{  T\otimes_{T}F(M)\ar[rd]_{\simeq}
\ar[rr]^{\Upsilon_{T,M}} & &  F(T\otimes_{T}M)\ar[ld]^{\simeq} \\
 & F(M) & } ,
\]
for every $M\in{}^{T}\mathcal{M}^{\coring{C}}.$ We have
$\Upsilon_{T,M}$ is natural in $T.$ By the theorem of Mitchell
\cite[Theorem 3.6.5]{Popescu:1973}, there exists a unique natural
transformation
\[
\Upsilon_{-,-}:-\otimes_{T}F(-)\rightarrow F(-\otimes_{T}-)
\]
extending the natural transformation $\Upsilon_{T,M}.$ We refer to
\cite{Gomez:2002} for more details.

\begin{remark}\label{Mitchell}
The theorem of Mitchell \cite[Theorem 3.6.5]{Popescu:1973} holds
also if we suppose only that the category $\mathscr{C'}$ is
preadditive and has coproducts, or if the category $\mathscr{C'}$
is preadditive and the functor $S$ preserves coproducts. This fact
is used to go up that the natural transformation $\Upsilon$ exists
for every $k$-linear functor
$F:\mathcal{M}^{\coring{C}}\rightarrow\mathcal{M}^{\coring{D}}$
even if the category $\mathcal{M}^{\coring{D}}$ is not abelian.
Note also that its corollary \cite[Corollary 3.6.6]{Popescu:1973}
holds also if we suppose only that the category $\mathscr{C'}$ is
preadditive.
\end{remark}

Let $M\in{}^{\coring{C}^{\prime}}\mathcal{M}^{\coring{C}}$ and
$N\in{}^{\coring{C}}\mathcal{M}^{\coring{C}^{\prime\prime}}$. The
map
\[
\omega_{M,N}=\rho_{M}\otimes_{A}N-M\otimes_{A}\lambda_{N}:M\otimes_{A}N
\rightarrow M\otimes_{A}\coring{C}\otimes_{A}N
\]
is a $\coring{C}^{\prime}-\coring{C}^{\prime\prime}$-bicomodule
map. Denote by $M\square_{\coring{C}}N$ its kernel in
$_{A^{\prime}}\mathcal{M}_{A^{\prime\prime}}$. If $\omega_{M,N}$
is $\coring{C}^{\prime}_{A^{\prime}}$-pure and
$_{A^{\prime\prime}}\coring{C}^{\prime\prime}$-pure, and the
following maps
\begin{equation}
\operatorname{ker}
(\omega_{M,N})\otimes_{A^{\prime\prime}}\coring{C}^{\prime\prime}
\otimes_{A^{\prime\prime}}\coring{C}^{\prime\prime}, \quad
\coring{C}^{\prime} \otimes_{A^{\prime}}\coring{C}^{\prime}
\otimes_{A^{\prime}}\operatorname{ker} (\omega_{M,N}) \quad
\textrm{and} \quad
\coring{C}^{\prime}\otimes_{A^{\prime}}\operatorname{ker}
(\omega_{M,N})\otimes_{A^{\prime\prime}}\coring{C}^{\prime\prime}
\end{equation}
are injective, then $M\square_{\coring{C}}N$ is the kernel of
$\omega_{M,N}$ in
${}^{\coring{C}^{\prime}}\mathcal{M}^{\coring{C}^{\prime\prime}}$.
This is the case if $\omega_{M,N}$ is
$(\coring{C}^{\prime}\otimes_{A^{\prime}}\coring{C}^{\prime})_{A^{\prime}}$-pure,
$_{A^{\prime\prime}}(\coring{C}^{\prime\prime}
\otimes_{A^{\prime\prime}}\coring{C}^{\prime\prime})$-pure, and
$\coring{C}^{\prime}\otimes_{A^{\prime}}\omega_{M,N}$ is
$_{A^{\prime\prime}}\coring{C}^{\prime\prime}$-pure (e.g. if
$\coring{C}^{\prime}_{A^{\prime}}$ and
$_{A^{\prime\prime}}\coring{C}^{\prime\prime}$  are flat, or if
$\coring{C}$ is a coseparable $A$-coring).

If for every
$M\in{}^{\coring{C}^{\prime}}\mathcal{M}^{\coring{C}}$ and
$N\in{}^{\coring{C}}\mathcal{M}^{\coring{C}^{\prime\prime}}$,
$\omega_{M,N}$ is $\coring{C}^{\prime}_{A^{\prime}}$-pure and
$_{A^{\prime\prime}}\coring{C}^{\prime\prime}$-pure, then we have
a bifunctor
\begin{equation}\label{cotensorbifunctor}
\xymatrix{-\cotensor{\coring{C}}-:{}^{\coring{C}^{\prime}}\mathcal{M}^{\coring{C}}
\times{}^{\coring{C}}\mathcal{M}^{\coring{C}^{\prime\prime}}\ar[r]&
{}^{\coring{C}^{\prime}}\mathcal{M}^{\coring{C}^{\prime\prime}}}.
\end{equation}
which is $k$-linear in each variable. If in particular
$\coring{C}^{\prime}_{A^{\prime}}$ and
$_{A^{\prime\prime}}\coring{C}^{\prime\prime}$ are flat, or if
$\coring{C}$ is a coseparable $A$-coring, then the bifunctor
\eqref{cotensorbifunctor} is well defined.

By a proof similar to that of \cite[II.1.3]{AlTakhman:1999}, we
have, for every $M\in\mathcal{M}^{\coring{C}}$, that the functor
$M\square_{\coring{C}}-$ preserves direct limits.

The following lemma was used implicitly in the proof of
\cite[Proposition 3.4]{Gomez:2002}, and it will be useful for us
in the proof of the next theorem.

\begin{lemma}\label{2}
If $M\in{}^{\coring{C}^{\prime}}\mathcal{M}^{\coring{C}}$, and
$F:\mathcal{M}^{\coring{C}}\rightarrow\mathcal{M}^{\coring{D}}$ is
a $M$-compatible $k$-linear functor in the sense of \cite[p.
210]{Gomez:2002}, which preserves
coproducts, then for all $X\in\mathcal{M}^{A^{\prime}},$%
\[
\Upsilon_{X\otimes_{A^{\prime}}\coring{C}^{\prime},M}(X\otimes_{A^{\prime}%
}\lambda_{F(M)})=F(X\otimes_{A^{\prime}}\lambda_{M})\Upsilon_{X,M}.
\]
\end{lemma}

\begin{proof}
Let us consider the diagram
\[
\xymatrix{
 X\otimes_{A^{\prime}}F(M) \ar[rrrr]^{X\otimes_{A^{\prime}}\lambda
_{F(M)}} \ar[rrd]^{X\otimes_{A^{\prime}}F(\lambda_{M})}
\ar[dd]_{\Upsilon_{X,M}}& & & & X
\otimes_{A^{\prime}}\coring{C}^{\prime }\otimes_{A^{\prime}}F(M)
\ar[lld]^{X\otimes_{A^{\prime}}\Upsilon_{\coring{C}^{\prime},M}}
\ar[dd]^{\Upsilon_{X\otimes
_{A^{\prime}}\coring{C}^{\prime},M}}\\
& & X\otimes_{A^{\prime}}F(\coring{C}
^{\prime}\otimes_{A^{\prime}}M)
\ar[rrd]^{\Upsilon_{X,\coring{C}^{\prime}\otimes_{A^{\prime}}M}} & & \\
F(X\otimes_{A^{\prime}}M)
\ar[rrrr]_{F(X\otimes_{A^{\prime}}\lambda _{M})}&  &  & &
F(X\otimes_{A^{\prime}}\coring{C}^{\prime }\otimes_{A^{\prime}}M)
}
\]
The commutativity of the top triangle follows from the definition
of $\lambda_{F(M)}$, while the right triangle commutes by
\cite[Lemma 3.3]{Gomez:2002} (we take
$S=T=A^{\prime}$, and $Y=\coring{C}^{\prime}%
$), and the left triangle is commutative since $\Upsilon_{X,-}$.
is natural. Therefore, the commutativity of the rectangle holds.
\end{proof}

 A closer analysis of \cite[Theorem 3.5]{Gomez:2002} gives the
following generalization of \cite[Proposition 2.1]{Takeuchi:1977}
and \cite[23.1(1)]{Brzezinski/Wisbauer:2003} (concerning this
last, the condition ``$F$ is kernel preserving'' is used in its
proof). Recall from \cite{Guzman:1989} that a coring $\coring{C}$
is said to be \emph{coseparable} if the comultiplication map
$\Delta_{\coring{C}}$ is a split monomorphism of
$\coring{C}$--bicomodules. Of course, the trivial $A$--coring
$\coring{C} = A$ is coseparable and, hencefoth, every result for
comodules over coseparable corings applies in particular for
modules over rings.

\begin{theorem}\label{3}
Let
$F:\mathcal{M}^{\coring{C}}\rightarrow\mathcal{M}^{\coring{D}}$ be
a $k$-linear functor, such that
\begin{enumerate}[(I)]
\item $_{B}\coring{D}$ is flat and $F$ preserves the kernel of $\rho_{N}\otimes_{A}\coring{C}-N\otimes_{A}%
\Delta_{\coring{C}}$ for every $N \in \rcomod{\coring{C}}$, or
\item $\coring{C}$ is a coseparable $A$-coring and the categories
$\mathcal{M}^{\coring{C}}$ and $\mathcal{M}^{\coring{D}}$ are
abelian.
\end{enumerate}
Assume that at least one of the following statements holds
\begin{enumerate}
\item $\coring{C}_{A}$ is projective, $F$ preserves coproducts,
and $\Upsilon_{N,\coring{C}}$, $\Upsilon_{N\otimes_{A}\coring{C}%
,\coring{C}}$ are isomorphisms for all
$N\in\mathcal{M}^{\coring{C}}$ (e.g. , if $A$ is semisimple and
$F$ preserves coproducts), or \item  $\coring{C}_{A}$ is flat, $F$
preserves direct limits, and $\Upsilon _{N,\coring{C}}$,
$\Upsilon_{N\otimes_{A}\coring{C},\coring{C}}$ are isomorphisms
for all $N\in\mathcal{M}^{\coring{C}}$ (e.g. , if $A$ is a von
Neumann regular ring and $F$ preserves direct limits), or \item
$F$ preserves inductive limits (e.g., if $F$ has a right adjoint).
\end{enumerate}
 Then $F$ is naturally equivalent to
$-\square_{\coring{C}}F(\coring{C}).$
\end{theorem}

\begin{proof}
At first, note that if $\coring{C}_{A}$ is projective, then the
right $A$-module $\coring{C}\otimes_{A}\coring{C}$ is projective
(by \cite[Example 3, p. 105, and Proposition
VI.9.5]{Stenstrom:1975}). Hence, if $\coring{C}_{A}$ is projective
and $F$ preserves coproducts, then $F$ is $M$-compatible in the
sense of \cite[p. 210]{Gomez:2002}, for all
$M\in{}^{\coring{C}}\mathcal{M}^{\coring{C}%
}.$ In each case, we have $F$ is $\coring{C}$-compatible where
$\coring{C}\in{}^{\coring{C}}\mathcal{M}^{\coring{C}}.$ Therefore,
by \cite[Proposition 3.4]{Gomez:2002}, $F(\coring{C})$ can be
viewed as a $\coring{C}-\coring{D}$-bicomodule. From Lemma
\ref{2}, and since $\Upsilon_{-,\coring{C}}$ is a natural
transformation, we have, for every $N\in
\mathcal{M}^{\coring{C}}$, the commutativity of the following
diagram with
exact rows in $\mathcal{M}^{\coring{D}}$%
\[%
\xymatrix{
 0 \ar[r] & N\square_{\coring{C}}F(\coring{C})\ar[r] &
 N\otimes_{A}F(\coring{C})
\ar[d]^{\simeq}_{\Upsilon_{N,\coring{C}}}
\ar[rrr]^-{\rho_{N}\otimes_{A}F(\coring{C})
-N\otimes_{A}%
\lambda_{F(\coring{C})}} &  & & N\otimes
_{A}\coring{C}\otimes_{A}F(\coring{C})\ar[d]^{\simeq
}_{\Upsilon_{N\otimes_{A}\coring{C},\coring{C}}} \\
 0 \ar[r] &  F(N)\ar[r]^{F(\rho
_{N})} & F(N\otimes_{A}\coring{C})
\ar[rrr]^-{F(\rho_{N}\otimes_{A}\coring{C}-N\otimes_{A}%
\Delta_{\coring{C}})} & & & F(N\otimes
_{A}\coring{C}\otimes_{A}\coring{C})}
\]
(the exactness of the bottom sequence for the case $(II)$ follows
by
factorizing the map $\omega_{N,\coring{C}} = \rho_{N}\otimes_{A}\coring{C}-N\otimes_{A}%
\Delta_{\coring{C}}$ through its image, and using the facts that
the sequence $\xymatrix{0 \ar[r] &  N\ar[r]^{\rho _{N}} & N
\otimes_{A}\coring{C} \ar[r]^-{\omega_{N,\coring{C}}}  & N \otimes
_{A}\coring{C}\otimes_{A}\coring{C}}$ is split exact in
$\rcomod{\coring{C}}$ in the sense of
\cite[40.5]{Brzezinski/Wisbauer:2003}, and that additive functors
between abelian categories preserve split exactness). By the
universal property of kernels, there exists a unique isomorphism
$\eta_{N}:N\square _{\coring{C}}F(\coring{C}) \rightarrow F(N)$ in
$\mathcal{M}^{\coring{D}}$ making commutative the above diagram.
It easy to show that $\eta$ is natural. Hence
$F\simeq-\square_{\coring{C}}F(\coring{C}).$
\end{proof}

As an immediate consequence of the last theorem we have the
following generalization of Eilenberg-Watts Theorem
\cite[Proposition VI.10.1]{Stenstrom:1975}.

\begin{corollary}\label{4}
Let
$F:\mathcal{M}^{\coring{C}}\rightarrow\mathcal{M}^{\coring{D}}$ be
a $k$-linear functor.\begin{enumerate}[1.]\item If
$_{B}\coring{D}$ is flat and $A$ is a semisimple ring (resp. a von
Neumann regular ring), then the following statements are
equivalent\begin{enumerate}[(a)]\item $F$ is left exact and
preserves coproducts (resp. left exact and preserves direct
limits); \item $F\simeq-\square_{\coring{C}}M$ for some bicomodule
$M\in{}^{\coring{C}%
}\mathcal{M}^{\coring{D}}$.\end{enumerate}
 \item If
$_{A}\coring{C}$ and $_{B}\coring{D}$ are flat, then the following
statements are equivalent \begin{enumerate}[(a)]\item $F$ is exact
and preserves inductive limits;\item
$F\simeq-\square_{\coring{C}}M$ for some
bicomodule $M\in{}^{\coring{C}%
}\mathcal{M}^{\coring{D}}$ which is coflat in
$^{\coring{C}}\mathcal{M}%
$. \end{enumerate} \item If $\coring{C}$ is a coseparable
$A$-coring and the categories $\mathcal{M}^{\coring{C}}$ and
$\mathcal{M}^{\coring{D}}$ are abelian, then the following
statements are equivalent \begin{enumerate}[(a)]\item $F$
preserves inductive
limits;\item $F$ preserves cokernels and $F\simeq-\square_{\coring{C}%
}M$ for some bicomodule
$M\in{}^{\coring{C}}\mathcal{M}^{\coring{D}}$.
\end{enumerate}\item If $\coring{C}=A$ and the category
$\mathcal{M}^{\coring{D}}$ is abelian, then the following
statements are equivalent
\begin{enumerate}[(a)]\item $F$ has a right adjoint;\item $F$
preserves inductive limits;\item
$F\simeq-\otimes_{A}M$ for some bicomodule $M\in{}^{A}\mathcal{M}%
^{\coring{D}}.$\end{enumerate}
\end{enumerate}
\end{corollary}

A bicomodule $N \in \bcomod{\coring{C}}{\coring{D}}$ is said to be
\emph{quasi-finite} as a right $\coring{D}$--comodule if the
functor $- \tensor{A} N : \rmod{A} \rightarrow
\rcomod{\coring{D}}$ has a left adjoint $h_{\coring{D}}(N,-) :
\rcomod{\coring{D}} \rightarrow \rmod{A}$, the \emph{co-hom
functor}. If $\omega_{Y,N}$ is $\coring{D} \tensor{B}
\coring{D}$--pure for every right $\coring{C}$--comodule $Y$
(e.g., ${}_B\coring{D}$ is flat or $\coring{C}$ is coseparable)
then $N_{\coring{D}}$ is quasi-finite if and only if $-
\cotensor{\coring{C}} N : \rcomod{\coring{C}} \rightarrow
\rcomod{\coring{D}}$ has a left adjoint, which we still to denote
by $h_{\coring{D}}(N,-)$ \cite[Proposition 4.2]{Gomez:2002}. The
particular case of the following statement when the co-hom is
exact generalizes \cite[Corollary 3.12]{AlTakhman:2002}.

\begin{corollary}\label{5}
Let $N\in{}^{\coring{C}}\mathcal{M}^{\coring{D}}$ be a bicomodule,
quasi-finite as a right $\coring{D}$-comodule, such that
$_{A}\coring{C}$ and $_{B}\coring{D}$ are flat. If the co-hom
functor $\cohom{\coring{D}}{N}{-}$ is exact or if $\coring{D}$ is
a coseparable $B$-coring, then we have
\[
\cohom{\coring{D}}{N}{-}
\simeq-\square_{\coring{D}}\cohom{\coring{D}}{N}{\coring{D}}:
\mathcal{M}^{\coring{D}}\rightarrow \mathcal{M}^{\coring{C}}.
\]
\end{corollary}

\begin{proof}
The functor $\cohom{\coring{D}}{N}{-}$ is $k$-linear and preserves
inductive limits, since it is a left adjoint to the $k$-linear
functor
$-\square_{\coring{C}}N:\mathcal{M}^{\coring{C}}\rightarrow
\mathcal{M}^{\coring{D}}$ (by \cite[Proposition 4.2]{Gomez:2002}).
Hence Theorem \ref{3} achieves the proof.
\end{proof}

Now we will use the following generalization of \cite[Lemma
2.2]{Takeuchi:1977}.

\begin{lemma}\label{18}
Let $\Lambda$, $\Lambda^{\prime}$ be
bicomodules in $^{\coring{D}%
}\mathcal{M}^{\coring{C}}$ and $G=-\square_{\coring{D}}\Lambda$,
$G^{\prime}=-\square_{\coring{D}}\Lambda^{\prime}.$ Suppose
moreover that ${}_A\coring{C}$ is flat and $B$ is a von Neumann
regular ring, or ${}_A\coring{C}$ is flat and $G$ and $G^{\prime}$
are cokernel preserving, or $\coring{D}$ is a coseparable coring.
Then
\[
\nat{G}{G^{\prime}}\simeq
\hom{\coring{D},\coring{C}}{\Lambda}{\Lambda^{\prime}}.
\]
\end{lemma}

\begin{proof}
Let $\alpha:G\rightarrow G^{\prime}$ be a natural transformation.
By \cite[Lemma 3.2(1)]{Gomez:2002}, $\alpha_{\coring{D}}$ is left
$B$-linear. For the rest of the proof it suffices to replace
$\otimes$ by $\otimes_{B}$ in the proof of \cite[Lemma
4.1]{Caenepeel/DeGroot/Militaru:2002}.
\end{proof}

The following proposition generalizes well-known facts from tensor
product functors for modules to cotensor product functors for
comodules over corings.

\begin{proposition}\label{19a}
Suppose that $_{A}\coring{C}$, $\coring{C}_{A}$, $_{B}\coring{D}$
and
$\coring{D}_{B}$ are flat. Let $X\in{}^{\coring{C}}\mathcal{M}%
^{\coring{D}}$ and
$\Lambda\in{}^{\coring{D}}\mathcal{M}^{\coring{C}}$. Consider the
following properties:
\begin{enumerate}[(1)]
\item \label{adj1} $-\square _{\coring{C}}X$ is left adjoint to
$-\square_{\coring{D}}\Lambda$; \item \label{adj2} $\Lambda$ is
quasi-finite as a right $\coring{C}$-comodule and $-
\square_{\coring{C}} X \simeq \cohom{\coring{C}}{\Lambda}{-}$;
\item \label{adjsep} $\Lambda$ is quasi-finite as a right $\coring{C}%
$-comodule and $X\simeq \cohom{\coring{C}}{\Lambda}{\coring{C}}$
in $^{\coring{C}}\mathcal{M}^{\coring{D}}$; \item
\label{semicontext} there exist bicolinear maps
\[
\psi:\coring{C}\rightarrow X\square_{\coring{D}}\Lambda\text{ and
}\omega:\Lambda\square_{\coring{C}}X\rightarrow\coring{D}
\]
in $^{\coring{C}}\mathcal{M}^{\coring{C}}$ and $^{\coring{D}}
\mathcal{M}^{\coring{D}}$ respectively, such that
\begin{equation}\label{unitcounit}
(\omega\square_{\coring{D}}\Lambda)\circ(\Lambda
\square_{\coring{C}}\psi)=\Lambda\text{ and }(X\square
_{\coring{D}}\omega)\circ(\psi\square_{\coring{C}}X)=X;
\end{equation}
\item \label{adj1symm} $\Lambda \cotensor{\coring{C}}-$ is left
adjoint to $X \cotensor{\coring{D}}-$.
\end{enumerate}
Then (\ref{adj1}) and (\ref{adj2}) are equivalent, and they imply
(\ref{adjsep}). The converse is true if $\coring{C}$ is a
coseparable $A$-coring. If ${}_A X$ and ${}_B \Lambda$ are flat,
and $\omega{}_{X,Y}=\rho{}_X\otimes_{B}Y-X\otimes_{A}\rho{}_Y$ is
pure as an $A$--linear map and
$\omega{}_{Y,X}=\rho{}_Y\otimes_{A}X-Y\otimes_{B}\rho{}_X$ is pure
as a $B$--linear map (e.g. if ${}_{\coring{C}}X$ and
${}_{\coring{D}}\Lambda$ are coflat
\cite[21.5]{Brzezinski/Wisbauer:2003} or $A$ and $B$ are von
Neumann regular rings), or if $\coring{C}$ and $\coring{D}$ are
coseparable, then (\ref{semicontext}) implies (\ref{adj1}). The
converse is true if ${}_{\coring{C}}X$ and
${}_{\coring{D}}\Lambda$ are coflat, or if $A$ and $B$ are von
Neumann regular rings, or if $\coring{C}$ and $\coring{D}$ are
coseparable. Finally, if $\coring{C}$ and $\coring{D}$ are
coseparable, or if $X$ and $\Lambda$ are coflat on both sides, or
if $A, B$ are von Neumann regular rings, then (\ref{adj1}),
(\ref{semicontext}) and (\ref{adj1symm}) are equivalent.
\end{proposition}

\begin{proof}
The equivalence between (\ref{adj1}) and (\ref{adj2}) follows from
\cite[Proposition 4.2]{Gomez:2002}. That (\ref{adj2}) implies
(\ref{adjsep}) is a consequence of \cite[Proposition
3.4]{Gomez:2002}. If $\coring{C}$ is coseparable and we assume
(\ref{adjsep}) then, by Corollary \ref{5},
$\cohom{\coring{C}}{\Lambda}{-} \simeq - \cotensor{\coring{C}}
\cohom{\coring{C}}{\Lambda}{\coring{C}} \simeq -
\cotensor{\coring{C}} X$. That (\ref{adj1}) implies
(\ref{semicontext}) follows from Lemma \ref{18} by evaluating the
unit and the counit of the adjunction at $\coring{C}$ and
$\coring{D}$, respectively. Conversely, if we put
$F=-\cotensor{\coring{C}}X$ and $G=-\cotensor{\coring{D}}\Lambda$,
we have $GF \simeq
-\square_{\coring{C}%
}(X\square_{\coring{D}}\Lambda)$ and $FG \simeq -\square
_{\coring{D}}(\Lambda\square_{\coring{C}}X)$ by \cite[Proposition
22.6]{Brzezinski/Wisbauer:2003}$. $ Define natural transformations
$$\xymatrix@1{\eta :1_{\mathcal{M}^{\coring{C}}}\ar[r]^{\simeq} &
-\square
_{\coring{C}}\coring{C}\ar[r]^-{-\square_{\coring{C}}\psi} & GF}$$
 and
$$ \xymatrix@1{\varepsilon:FG\ar[r]^{-\square_{\coring{D}}\omega} &
-\square_{\coring{D}}\coring{D}\ar[r]^{\simeq} &
1_{\mathcal{M}^{\coring{D}}},} $$ which become the unit and the
counit of an adjunction by \eqref{unitcounit}. This gives the
equivalence between (\ref{adj1}) and (\ref{semicontext}). The
equivalence between (\ref{semicontext}) and (\ref{adj1symm})
follows by symmetry.
\end{proof}

\begin{definition}\label{7}
Following \cite{AlTakhman:2002} and
\cite{Brzezinski/Wisbauer:2003}, a bicomodule
$N\in{}^{\coring{C}}\mathcal{M}^{\coring{D}}$ is called an
\emph{injector} (resp. \emph{an injector-cogenerator}) as a right
$\coring{D}$-comodule if the functor
$-\otimes_{A}N:\mathcal{M}^{A}\rightarrow\mathcal{M}^{\coring{D}}$
preserves injective (resp. injective cogenerator) objects.
\end{definition}

\begin{proposition}\label{19}
Suppose that $_{A}\coring{C}$ and $_{B}\coring{D}$ are flat. Let
$X\in{}^{\coring{C}}\mathcal{M}%
^{\coring{D}}$ and
$\Lambda\in{}^{\coring{D}}\mathcal{M}^{\coring{C}}.$ The following
statements are equivalent
\begin{enumerate}[(i)]
\item $-\square _{\coring{C}}X$ is left adjoint to
$-\square_{\coring{D}}\Lambda$, and $- \cotensor{\coring{C}} X$ is
left exact (or ${}_AX$ is flat or ${}_{\coring{C}}X$ is coflat);
\item $\Lambda$ is quasi-finite as a right $\coring{C}$-comodule,
$- \square_{\coring{C}} X \simeq \cohom{\coring{C}}{\Lambda}{-}$,
and $- \square_{\coring{C}}X$ is left exact (or ${}_AX$ is flat or
${}_{\coring{C}}X$ is coflat); \item
 $\Lambda$ is quasi-finite and injector as a right
$\coring{C}$-comodule and $X\simeq
\cohom{\coring{C}}{\Lambda}{\coring{C}}$ in
$^{\coring{C}}\mathcal{M}^{\coring{D}}$.
\end{enumerate}
\end{proposition}

\begin{proof}
First, observe that if ${}_{\coring{C}}X$ is coflat, then ${}_AX$
is flat \cite[21.6]{Brzezinski/Wisbauer:2003}, and that if ${}_AX$
is flat, then the functor $- \cotensor{\coring{C}} X$ is left
exact. Thus, in view of Proposition \ref{19a}, it suffices if we
prove that the version of $(ii)$ with $- \cotensor{\coring{C}} X$
left exact implies $(iii)$, and this last implies the version of
$(ii)$ with ${}_{\coring{C}}X$ coflat. Assume that $-
\cotensor{\coring{C}}X \simeq \cohom{\coring{C}}{\Lambda}{-}$ with
$- \cotensor{\coring{C}}X$ left exact. By \cite[Proposition
3.4]{Gomez:2002}, $X \simeq
\cohom{\coring{C}}{\Lambda}{\coring{C}}$ in
$\bcomod{\coring{C}}{\coring{D}}$. Being a left adjoint,
$\cohom{\coring{C}}{\Lambda}{-}$ is right exact and, henceforth,
exact. By \cite[Theorem 3.2.8]{Popescu:1973},
$\Lambda_{\coring{C}}$ is an injector and we have proved $(iii)$.
Conversely, if $\Lambda_{\coring{C}}$ is a quasi-finite injector
and $X \simeq \cohom{\coring{C}}{\Lambda}{\coring{C}}$ as
bicomodules, then $ - \cotensor{\coring{C}} X \simeq -
\cotensor{\coring{C}} \cohom{\coring{C}}{\Lambda}{\coring{C}}$
and, by \cite[Theorem 3.2.8]{Popescu:1973}, we get that
$\cohom{\coring{C}}{\Lambda}{-}$ is an exact functor. By Corollary
\ref{5}, $\cohom{\coring{C}}{\Lambda}{-} \simeq -
\cotensor{\coring{C}} \cohom{\coring{C}}{\Lambda}{\coring{C}}
\simeq - \cotensor{\coring{C}} X$, and ${}_{\coring{C}}X$ is
coflat.
\end{proof}

We are ready to give our characterization of Frobenius functors
between categories of comodules over corings. We left to the
reader to derive the corresponding characterizations for the
particular cases of rings or of coalgebras over fields.

\begin{theorem}\label{20a}
Suppose that $_{A}\coring{C}$ and $_{B}\coring{D}$ are flat. Let
$X\in{}^{\coring{C}}\mathcal{M}%
^{\coring{D}}$ and
$\Lambda\in{}^{\coring{D}}\mathcal{M}^{\coring{C}}.$ The following
statements are equivalent
\begin{enumerate}[(i)]
\item $(-\cotensor{\coring{C}} X, - \cotensor{\coring{D}}
\Lambda)$ is a Frobenius pair; \item $ - \cotensor{\coring{C}} X$
is a Frobenius functor, and $\cohom{\coring{D}}{X}{\coring{D}}
\simeq \Lambda$ as bicomodules; \item there is a Frobenius pair
$(F,G)$ for $\rcomod{\coring{C}}$ and $\rcomod{\coring{D}}$ such
that $F(\coring{C}) \simeq \Lambda$ and $G(\coring{D}) \simeq X$
as bicomodules; \item $\Lambda_{\coring{C}}, X_{\coring{D}}$ are
quasi-finite injectors,  and $X \simeq
\cohom{\coring{C}}{\Lambda}{\coring{C}}$ and $\Lambda \simeq
\cohom{\coring{D}}{X}{\coring{D}}$ as bicomodules; \item
$\Lambda_{\coring{C}}, X_{\coring{D}}$ are quasi-finite, and $-
\cotensor{\coring{C}} X \simeq \cohom{\coring{C}}{\Lambda}{-}$ and
$ - \cotensor{\coring{D}} \Lambda \simeq
\cohom{\coring{D}}{X}{-}$.
\end{enumerate}
\end{theorem}

\begin{proof}
$(i) \Leftrightarrow (ii) \Leftrightarrow (iii)$ This is obvious,
after Theorem \ref{3} and \cite[Proposition 3.4
]{Gomez:2002}. \\
$(i) \Leftrightarrow (iv)$ Follows from Proposition \ref{19}.\\
$(iv) \Leftrightarrow (v)$. If $X_{\coring{D}}$ and
$\Lambda_{\coring{C}}$ are quasi-finite, then ${}_AX$ and
${}_B\Lambda$ are flat. Now, apply Proposition \ref{19}.
\end{proof}

From Proposition \ref{19a} and Proposition \ref{19} (or Theorem
\ref{20a}) we get the following

\begin{theorem}\label{leftrightFrobenius}
Let $X\in{}^{\coring{C}}\mathcal{M}%
^{\coring{D}}$ and
$\Lambda\in{}^{\coring{D}}\mathcal{M}^{\coring{C}}.$ Suppose that
$_{A}\coring{C}$, $\coring{C}_{A}$, $_{B}\coring{D}$ and
$\coring{D}_{B}$ are flat. The following statements are equivalent
\begin{enumerate} \item $(-\cotensor{\coring{C}} X, -
\cotensor{\coring{D}}\Lambda)$ is a Frobenius pair, with
$X_{\coring{D}}$ and $\Lambda_{\coring{C}}$ coflat; \item
$(\Lambda \cotensor{\coring{C}} -, X \cotensor{\coring{D}}-)$ is a
Frobenius pair, with ${}_{\coring{C}}X$ and
${}_{\coring{D}}\Lambda$ coflat; \item $X$ and $\Lambda$ are
coflat quasi-finite injectors on both sides, and $X\simeq
\cohom{\coring{C}}{\Lambda}{\coring{C}}$ in
$^{\coring{C}}\mathcal{M}^{\coring{D}}$ and $\Lambda\simeq
\cohom{\coring{D}}{X}{\coring{D}}$ in
$^{\coring{D}}\mathcal{M}^{\coring{C}}$.
 \end{enumerate}
  If moreover $\coring{C}$ and $\coring{D}$ are coseparable (resp. $A$
and $B$ are von Neumann regular rings), then the following
statements are equivalent
\begin{enumerate}
\item $(-\cotensor{\coring{C}} X,- \cotensor{\coring{D}}\Lambda)$
is a Frobenius pair; \item $(\Lambda \cotensor{\coring{C}} -, X
\cotensor{\coring{D}}-)$ is a Frobenius pair;\item $X$ and
$\Lambda$ are quasi-finite (resp. quasi-finite injector) on both
sides, and $X\simeq \cohom{\coring{C}}{\Lambda}{\coring{C}}$ in
$^{\coring{C}}\mathcal{M}^{\coring{D}}$ and $\Lambda\simeq
\cohom{\coring{D}}{X}{\coring{D}}$ in
$^{\coring{D}}\mathcal{M}^{\coring{C}}$.
\end{enumerate}
\end{theorem}

\section{Frobenius functors between corings with a duality}

We will look to Frobenius functors for corings closer to
coalgebras over fields, in the sense that the categories of
comodules share a fundamental duality.

\par

An object $M$ of a Grothendieck category $\cat{C}$ is said to be
\emph{finitely generated} \cite[p. 121]{Stenstrom:1975} if
whenever $M = \sum_i M_i$ is a direct union of subobjects $M_i$,
then $M = M_{i_0}$ for some index $i_0$. Alternatively, $M$ is
finitely generated if the functor $\hom{\cat{C}}{M}{-}$ preserves
direct unions \cite[Proposition V.3.2]{Stenstrom:1975}. The
category $\cat{C}$ is \emph{locally finitely generated} if it has
a family of finitely generated generators. Recall from \cite[p.
122]{Stenstrom:1975} that a finitely generated object $M$ is
\emph{finitely presented} if every epimorphism $L \rightarrow M$
with $L$ finitely generated has finitely generated kernel. By
\cite[Proposition V.3.4]{Stenstrom:1975}, if $\cat{C}$ is locally
finitely generated, then $M$ is finitely presented if and only if
$\hom{\cat{C}}{M}{-}$ preserves direct limits. For the notion of a
locally projective module we refer to
\cite{Zimmermann-Huisgen:1976}.

\begin{lemma}\label{locally}
Let $\coring{C}$ be a coring over a ring $A$ such that
${}_A\coring{C}$ is flat.
\begin{enumerate}[(1)]
\item\label{fg} A comodule $M \in \rcomod{\coring{C}}$ is finitely
generated if and only if $M_A$ is finitely generated.
\item\label{fp} A comodule $M \in \rcomod{\coring{C}}$ is finitely
presented if $M_A$ is finitely presented. The converse is true
whenever $\rcomod{\coring{C}}$ is locally finitely generated.
\item\label{lfg} If ${}_A\coring{C}$ is locally projective, then
$\rcomod{\coring{C}}$ is locally finitely generated.
\end{enumerate}
\end{lemma}

\begin{proof}
The forgetful functor $U : \rcomod{\coring{C}} \rightarrow
\rmod{A}$  has an exact left adjoint $ - \tensor{A} \coring{C} :
\rmod{A} \rightarrow \rcomod{\coring{C}}$ which preserves direct
limits. Thus, $U$ preserves finitely generated objects and, in
case that $\rcomod{\coring{C}}$ is locally finitely generated,
finitely presented objects. Now, if $M \in \rcomod{\coring{C}}$ is
finitely generated as a right $A$-module, and $M = \sum_i M_i$ as
a direct union of subcomodules, then $U(M) = U(\sum_iMi) =
\sum_iU(M_i)$, since $U$ is exact and preserves coproducts.
Therefore, $U(M) = U(M_{i_0})$ for some index $i_0$ which implies,
being $U$ a faithfully exact functor, that $M = M_{i_0}$. Thus,
$M$ is a finitely generated comodule. We have thus proved
(\ref{fg}), and the converse to (\ref{fp}). Now, if $M \in
\rcomod{\coring{C}}$ is such that $M_A$ is finitely presented,
then for every exact sequence $0 \rightarrow K \rightarrow L
\rightarrow M \rightarrow 0$ in $\rcomod{\coring{C}}$ with $L$
finitely generated, we get an exact sequence $0 \rightarrow K_A
\rightarrow L_A \rightarrow M_A \rightarrow 0$ with $M_A$ finitely
presented. Thus, $K_A$ is finitely generated and, by (\ref{fg}),
$K \in \rcomod{\coring{C}}$ is finitely generated. This proves
that $M$ is a finitely presented comodule. Finally, (\ref{lfg}) is
a consequence of (\ref{fg}) and
\cite[19.12(1)]{Brzezinski/Wisbauer:2003}.
\end{proof}

The notation $\cat{C}_f$ stands for the full subcategory of a
Grothendieck category $\cat{C}$ whose objects are the finitely
generated objects. The category $\cat{C}$ is locally noetherian
\cite[p. 123]{Stenstrom:1975} if it has a family of noetherian
generators or equivalently, if $\mathbf{C}$ is locally finitely
generated and every finitely generated object of $\mathbf{C}$ is
noetherian. By \cite[Propositon V.4.2, Proposition V.4.1, Lemma
V.3.1(i)]{Stenstrom:1975}, in an arbitrary Grothendieck category,
every finitely generated object is noetherian if and only if every
finitely generated object is finitely presented. The version for
categories of modules of the following result is well-known.

\begin{lemma}\label{10}
Let $\cat{C}$ be a locally finitely generated category.
\begin{enumerate}[(1)]
\item\label{additive} The category $\cat{C}_f$ is additive.
\item\label{cokernels} The category $\cat{C}_f$ has cokernels, and
every monomorphism in $\cat{C}_f$ is a monomorphism in $\cat{C}$.
\item\label{locnoeth} The following statements are equivalent:
  \begin{enumerate}[(a)]
   \item\label{haskernels} The category $\cat{C}_f$ has kernels;
   \item\label{locnoeth2} $\cat{C}$ is locally noetherian;
   \item\label{abelian} $\cat{C}_f$ is abelian;
   \item\label{abeliansub} $\cat{C}_f$ is an abelian subcategory of
$\cat{C}$.
  \end{enumerate}
\end{enumerate}
\end{lemma}

\begin{proof}
1. Straightforward.\newline 2. That $\mathbf{C}_{f}$ has cokernels
is straightforward from \cite[Lemma V.3.1(i)]{Stenstrom:1975}.
Now, let $f:M\rightarrow N$ be a monomorphism in $\mathbf{C}_{f}$
and $\xi:X\rightarrow M$ be a morphism in $\mathbf{C}$ such that
$f\xi=0.$ Suppose that $X=\bigcup_{i\in I}X_{i}$, where
$X_{i}\in\mathbf{C}_{f}$, and $\iota_{i}:X_{i}\rightarrow X$,
$i\in I$ the canonical injections$.$ Then $f\xi\iota_{i}=0$, and
$\xi\iota_{i}=0$,
for every $i$, and by the definition of the inductive limit, $\xi=0.$%
\newline 3. $(b)\Rightarrow(a)$
Straightforward from \cite[Proposition
V.4.1]{Stenstrom:1975}.\newline $(d)\Rightarrow(c)$ and
$(c)\Rightarrow(a)$ are trivial.\newline $(a)\Rightarrow(b)$ Let
$M\in\mathbf{C}_{f}$, and $K$ be a subobject of $M.$ Let $\iota
:L\rightarrow M$ the kernel of the canonical morphism
$f:M\rightarrow M/K$ in $\mathbf{C}_{f}.$ Suppose that
$K=\bigcup_{i\in I}K_{i}$, where $K_{i}\in\mathbf{C}_{f}$, for
every $i\in I.$ By the universal property of the kernel, there
exist a unique morphism $\alpha:L\rightarrow K,$ and a unique
morphism $\beta_{i}:K_{i}\rightarrow L$, for every $i\in I$,
making commutative the diagrams
\[%
\xymatrix{
&   K_{i} \ar[d] \ar[dl]_{\beta_{i}} & \\
L \ar[dr]_{\alpha} \ar[r]^{\iota} & M \ar[r]^{f} & M/K \\
   & K \ar[u] &   }
\]
By (2), $\iota$ is a monomorphism in $\mathbf{C}$, then for every
$K_{i}\subset K_{j}$, the diagram
\[%
\xymatrix{
K_{i}  \ar[d] \ar[r]^{\beta_{i}} & L\\
\ar[ur]_{\beta_{j}}  K_{j} &   }
\]
commutes. Therefore we have the commutative diagram
\[%
\xymatrix{
& K \ar[d] \ar[dl]_{{\lim \atop {\longrightarrow \atop I}}\beta_{i}}\\
L \ar[r]^{\iota} & M %
}
\]
Then $K\simeq L$, and hence $K\in\mathbf{C}_{f}.$ Finally, by
\cite[Proposition V.4.1]{Stenstrom:1975}, $M$ is noetherian in
$\mathbf{C.}$\newline $(b)\Rightarrow(d)$ Straightforward from
\cite[Theorem 3.41]{Freyd:1964}.
\end{proof}

The following is a generalization of \cite[Proposition
3.1]{Castano/Gomez/Nastasescu:1999}.

\begin{proposition}\label{11}
Let $\mathbf{C}$ and $\mathbf{D}$ be two locally noetherian
categories. Then
\begin{enumerate}[(1)]
\item If $F:\mathbf{C}\rightarrow\mathbf{D}$ is a Frobenius
functor, then its restriction $F_{f}:\mathbf{C}_{f}\rightarrow
\mathbf{D}_{f}$ is a Frobenius functor.
\item If $H:\mathbf{C}%
_{f}\rightarrow\mathbf{D}_{f}$ is a Frobenius functor, then $H$
can be
uniquely extended to a Frobenius functor $\overline{H}:\mathbf{C}%
\rightarrow\mathbf{D}.$ \item The assignment $F\mapsto F_{f}$
defines a bijective correspondence (up to natural isomorphisms)
between Frobenius functors from $\mathbf{C}$ to $\mathbf{D}$ and
Frobenius functors from $\mathbf{C}_{f}$ to $\mathbf{D}_{f}.$
\item In particular, if
$\mathbf{C}=\mathcal{M}^{\coring{C}}$ and $\mathbf{D}=\mathcal{M}%
^{\coring{D}}$ are locally noetherian such that $_{A}\coring{C}$
and
$_{B}\coring{D}$ are flat, then $F:\mathcal{M}^{\coring{C}%
}\rightarrow\mathcal{M}^{\coring{D}}$ is a Frobenius functor if
and only if it preserves direct limits and comodules which are
finitely generated as right $A$-modules,
and the restriction functor $F_{f}:\mathcal{M}_{f}^{\coring{C}}%
\rightarrow\mathcal{M}_{f}^{\coring{D}}$ is a Frobenius functor.
\end{enumerate}
\end{proposition}

\begin{proof}
The proofs of \cite[Proposition 3.1 and Remark
3.2]{Castano/Gomez/Nastasescu:1999} remain valid for our
situation, but with some minor modifications: to prove that
$\overline{H}$ is well-defined, we use Lemma \ref{10}. (In the
proof of the statements (1), (2) and (3) we use the Grothendieck
AB-5 condition).
\end{proof}

In order to generalize \cite[Proposition A.2.1]{Takeuchi:1977ca}
and its proof, we need

\begin{lemma}\label{12}
\begin{enumerate}[(1)]
\item Let $\mathbf{C}$ be a locally noetherian category, let
$\mathbf{D}$ be an arbitrary Grothendieck category,
$F:\mathbf{C}\rightarrow\mathbf{D}$ be an
arbitrary functor which preserves direct limits, and $F_{f}:\mathbf{C}%
_{f}\rightarrow\mathbf{D}$ be its restriction to $\mathbf{C}_{f}.$
Then $F$ is exact (faithfully exact, resp. left, right exact) if
and only if $F_{f}$ is exact (faithfully exact, resp. left, right
exact). \\ In particular, an object $M$ in $\cat{C}_f$ is
projective (resp. projective generator) if and only if it is
projective (resp. projective generator) in $\cat{C}$. \item Let
$\mathbf{C}$ be a locally noetherian category. For every object
$M$ of $\mathbf{C}$, the following conditions are equivalent
\begin{enumerate}[(a)]
\item $M$ is injective (resp. an injective cogenerator); \item the
contravariant functor $\hom{\mathbf{C}}{-}{M}
:\mathbf{C}\rightarrow \mathbf{Ab}$ is exact (resp. faithfully
exact); \item the contravariant
functor $\hom{\mathbf{C}}{-}{M}_{f}:\mathbf{C}_{f}%
\rightarrow\mathbf{Ab}$ is exact (resp. faithfully exact).\\
In particular, an object $M$ in $\cat{C}_f$ is injective (resp.
injective cogenerator) if and only if it is injective (resp.
injective cogenerator) in $\cat{C}$.
\end{enumerate}
\end{enumerate}
\end{lemma}

\begin{proof}
(1) The ``only if'' part is straightforward from the fact that the
injection functor $\mathbf{C}_{f}\rightarrow\mathbf{C}$ is
faithfully exact.

 For the ``if'' part, suppose
that $F_{f}$ is left exact. Let $f:M\rightarrow N$ be a morphism
in $\mathbf{C}$. Put $M=\bigcup_{i\in I}M_{i}$ and
$N=\bigcup_{j\in J}N_{j},$ as direct union of directed families of
finitely generated subobjects. For $(i,j)\in I\times J,$ let $M_{i,j}%
=M_{i}\cap f^{-1}(N_{j})$, and $f_{i,j}:M_{i,j}\rightarrow N_{j}$
be the restriction of $f$ to $M_{i,j}.$ We have
$f=\underset{\underset{I\times J}{\longrightarrow}}{\lim}f_{i,j}$
and then $F(f)=\underset {\underset{I\times
J}{\longrightarrow}}{\lim}F_{f}(f_{i,j}).$ Hence\newline $\ker
F(f)=\ker\underset{\underset{I\times
J}{\longrightarrow}}{\lim}F_{f}(f_{i,j})=\underset
{\underset{I\times J}{\longrightarrow}}{\lim}\ker F_{f}( f_{i,j})
=\underset{\underset{I\times J}{\longrightarrow}}{\lim}F_{f}(\ker
f_{i,j})=\underset{\underset{I\times
J}{\longrightarrow}}{\lim}F(\ker f_{i,j})$ (by Lemma \ref{10})
$=F(\underset{\underset{I\times
J}%
{\longrightarrow}}{\lim}\ker f_{i,j})=F(\ker f).$\newline Finally
$F$ is left exact. Analogously, it can be proved that $F_{f}$ is
right exact implies that $F$ is also right exact. Now, suppose
that $F_{f}$ is faithfully exact. We have already proved that $F$
is exact. It remains to prove that $F$\ is faithful. For this, let
$0\neq M=\bigcup_{i\in I}M_{i}$ be an object of $\mathbf{C}$,
where $M_{i}$ is finitely generated for every $i\in I.$ We have
\[
F(M)=\underset{\underset{I}{\longrightarrow}}{\lim}F_{f}(M_{i})
\simeq\sum_{I}F_{f}(M_{i})
\]
(since $F$ is exact). Since $M\neq0$, there exists some $i_{0}\in
I$ such that $M_{i_{0}}\neq0.$ By \cite[Proposition
IV.6.1]{Stenstrom:1975}, $F_{f}(M_{i_{0}})\neq0,$ hence
$F(M)\neq0.$ Also by \cite[Proposition IV.6.1]{Stenstrom:1975},
$F$ is faithful.
\newline (2). $(a)\Leftrightarrow(b)$ Obvious.
\newline $(b)\Rightarrow(c)$
Analogous to that of the ``only if'' part of (1).\newline
$(c)\Rightarrow(a)$ Suppose that the functor
$\hom{\mathbf{C}}{-}{M}_{f}$ is exact. Since $\mathbf{C}$ is
locally noetherian, and by \cite[Proposition V.2.9 and Proposition
V.4.1]{Stenstrom:1975}, $M$ is injective. Now, suppose moreover
that $\hom{\mathbf{C}}{-}{M}_{f}$ is faithful. Let $L$ be a
non-zero object of $\mathbf{C}$, and $K$ be a non-zero finitely
generated subobject of $L.$ By \cite[Proposition
IV.6.1]{Stenstrom:1975}, there exists a non-zero morphism
$K\rightarrow M$. Since $M$ is injective, there exists a non-zero
morphism $L\rightarrow M$ making commutative the following diagram
\[%
\xymatrix{
&   M   & \\
0  \ar[r] & K \ar[r] \ar[u] & L. \ar[ul]}
\]
From \cite[Proposition IV.6.5]{Stenstrom:1975}, it follows that
$M$\ is a cogenerator.
\end{proof}

If $\coring{C}_A$ is flat and $M \in \rcomod{\coring{C}}$ is
finitely presented as right $A$-module, then
\cite[19.19]{Brzezinski/Wisbauer:2003} the dual left $A$-module
$\rdual{M} = \hom{A}{M}{A}$ has a left $\coring{C}$-comodule
structure
$$\rdual{M} \simeq
\hom{\coring{C}}{M}{\coring{C}} \subseteq \hom{A}{M}{\coring{C}}
\simeq \coring{C} \tensor{A} \rdual{M}.$$ Now, if ${}_A\rdual{M}$
turns out to be finitely presented and ${}_A\coring{C}$ is flat,
then $\ldual{(\rdual{M})} = \hom{A}{\rdual{M}}{A}$ is a right
$\coring{C}$-comodule and the canonical map $\sigma_M : M
\rightarrow \ldual{(\rdual{M})}$ is a homomorphism in
$\rcomod{\coring{C}}$. This construction leads to a duality (i.e.
a contravariant equivalence)
\[
\rDual: \rcomod{\coring{C}}_0 \leftrightarrows
\lcomod{\coring{C}}_0: \lDual
\]
between the full subcategories $\rcomod{\coring{C}}_0$ and
$\lcomod{\coring{C}}_0$ of $\rcomod{\coring{C}}$ and
$\lcomod{\coring{C}}$ whose objects are the comodules which are
finitely generated and projective over $A$ on the corresponding
side (this holds even without flatness assumptions of
$\coring{C}$). Call it \emph{basic duality} (details may be found
in \cite{Caenepeel/DeGroot/Vercruysse:unp}). Of course, in the
case $A$ is semisimple (e.g. for coalgebras over fields) these
categories are that of finitely generated comodules, and this
basic duality plays a remarkable role in the study of several
notions in the coalgebra setting (e.g. Morita equivalence
\cite{Takeuchi:1977}, semiperfect coalgebras \cite{Lin:1977},
Morita duality \cite{Gomez/Nastasescu:1995},
\cite{Gomez/Nastasescu:1996}, or Frobenius Functors
\cite{Castano/Gomez/Nastasescu:1999}). It would be interesting to
know, in the coring setting, to what extent the basic duality can
be extended to the subcategories $\rcomod{\coring{C}}_f$ and
$\lcomod{\coring{C}}_f$, since, as we will try to show in this
section, this allows to obtain better results. Of course, this is
the underlying idea when the ground ring $A$ is assumed to be
Quasi-Frobenius (see e.g. \cite{ElKaoutit/Gomez:2002unp}) for the
case of semiperfect corings and Morita duality), but we hope
future developments of the theory would have some pay off from the
more general setting we propose here.

Consider contravariant functors between Grothendieck categories $
H: \cat{A} \leftrightarrows \cat{A}': H', $ together with natural
transformations $\tau : 1_{\cat{A}} \rightarrow H' \circ H$ and
$\tau' : 1_{\cat{A}'} \rightarrow H \circ H'$, satisfying the
condition $H(\tau_X) \circ \tau'_{H(X)} = 1_{H(X)}$ and
$H'(\tau'_{X'}) \circ \tau_{H'(X')} = 1_{H'(X')}$ for $X \in
\cat{A}$ and $X' \in \cat{A}'$. Following
\cite{Colby/Fuller:1983}, this situation is called a \emph{right
adjoint pair}.

\begin{proposition}\label{lnoethrightadj}
Let $\coring{C}$ be an $A$-coring such that ${}_A\coring{C}$ and
$\coring{C}_A$ are flat. Assume that $\rcomod{\coring{C}}$ and
$\lcomod{\coring{C}}$ are locally noetherian categories. If
${}_A\rdual{M}$ and $\ldual{N}_A$ are finitely generated modules
for every $M \in \rcomod{\coring{C}}_f$ and $N \in
\lcomod{\coring{C}}_f$, then the basic duality extends to a right
adjoint pair $ \rDual: \rcomod{\coring{C}}_f \leftrightarrows
\lcomod{\coring{C}}_f: \lDual$.
\end{proposition}

\begin{proof}
If $M \in \rcomod{\coring{C}}_f$ then, since $\rcomod{\coring{C}}$
is locally noetherian, $M_{\coring{C}}$ is finitely presented. By
Lemma \ref{locally}, $M_A$ is finitely presented and the left
$\coring{C}$-comodule $\rdual{M}$ makes sense. Now, the assumption
${}_A\rdual{M}$ finitely generated implies, by Lemma
\ref{locally}, that $\rdual{M} \in \lcomod{\coring{C}}_f$. We have
then the functor $\rDual : \rcomod{\coring{C}}_f \rightarrow
\lcomod{\coring{C}}_f$. The functor $\rDual$ is analogously
defined, and the rest of the proof consists of straightforward
verifications (see \cite[Proposition
20.14(1)]{Anderson/Fuller:1992}, where $\tau$ and  $\tau'$ are the
evaluation maps).
\end{proof}

\begin{example}
The hypotheses are fulfilled if ${}_A\coring{C}$ and
$\coring{C}_A$ are locally projective and $A$ is left and right
noetherian (in this case the right adjoint pair already appears in
\cite{ElKaoutit/Gomez:2002unp}). But there are situations in which
no finiteness condition need to be required to $A$: this is the
case, for instance, of cosemisimple corings (see \cite[Theorem
3.1]{ElKaoutit/Gomez/Lobillo:2001unp}). In particular, if an
arbitrary ring $A$ contains a division ring $B$, then, by
\cite[Theorem 3.1]{ElKaoutit/Gomez/Lobillo:2001unp} the canonical
coring $A \tensor{B} A$ satisfies all hypotheses in Proposition
\ref{lnoethrightadj}.
\end{example}

\begin{definition}
Let $\coring{C}$ be a coring over $A$ satisfying the assumptions
of Proposition \ref{lnoethrightadj}. We will say that $\coring{C}$
\emph{has a duality} if the basic duality extends to a duality $$
\rDual: \rcomod{\coring{C}}_f \leftrightarrows
\lcomod{\coring{C}}_f: \lDual$$
\end{definition}

We have the following examples of a coring which has a
duality:\begin{itemize} \item $\coring{C}$ is a coring over a QF
ring $A$ such that ${}_A\coring{C}$ and $\coring{C}_A$ are flat
(and hence projective);\item $\coring{C}$ is a cosemisimple
coring;\item $\coring{C}$ is a coring over $A$ such that
${}_A\coring{C}$ and $\coring{C}_A$ are flat and semisimple,
$\rcomod{\coring{C}}$ and $\lcomod{\coring{C}}$ are locally
noetherian categories, and the dual of every simple right (resp.
left) $A$-module in the decomposition of $\coring{C}_A$ (resp.
${}_A\coring{C}$) as a direct sum of simple $A$-modules is
finitely generated and $A$-reflexive (in fact, every right (resp.
left) $\coring{C}$-comodule $M$ is a submodule of the semisimple
right (resp.left) $A$-module $M\otimes_{A}\coring{C}$, and hence
$M_A$ (resp. ${}_AM$) is also semisimple.)\end{itemize}

\begin{proposition}\label{13}
Suppose that the coring $\coring{C}$ has a duality. Let
$M\in\mathcal{M}^{\coring{C}}$ such that $M_{A}$ is flat. The
following are equivalent \begin{enumerate} \item $M$ is coflat
(resp. faithfully coflat); \item
$\hom{\coring{C}}{-}{M}_{f}:\mathcal{M}_{f}^{\coring{C}}\rightarrow\mathcal{M}_{k}$
is exact (resp. faithfully exact); \item $M$ is injective (resp.
an injective cogenerator).\end{enumerate}
\end{proposition}

\begin{proof}
Let $M\in\mathcal{M}%
^{\coring{C}}$ and $N\in{}^{\coring{C}}\mathcal{M}_{f}.$ We have
the following commutative diagram (in $\mathcal{M}_{k}$)
\[%
\xymatrix{  0 \ar[r] & M\square_{\coring{C}}N  \ar[r] &
M\otimes_{A}N \ar[rrr]^{\rho_{M}\otimes_{A}N-M\otimes_{A}\lambda_{N}%
} \ar[d] & & & M\otimes_{A}\coring{C}\otimes_{A}N \ar[d] \\
 0 \ar[r] & \hom{\coring{C}}{N^{\ast}}{M}\ar[r] &
\hom{A}{N^{\ast}}{M}\ar[rrr]^{f\mapsto\rho_{M}f-(f\otimes_{A}\coring{C})\rho_{N^{\ast}}
} & & & \hom{A}{N^{\ast}}{M\otimes_{A}\coring{C}}}
\]
where the vertical maps are the canonical maps. By the universal
property of the kernel, there is a unique morphism
$\eta_{M,N}:M\square_{\coring{C}}N\rightarrow
\hom{\coring{C}}{N^{\ast}}{M}$ making commutative the above
diagram. By the cube Lemma (see \cite[p. 43]{Mitchell:1965}),
$\eta$ is a natural transformation of bifunctors. If $M_{A}$ is
flat then $\eta_{M,N}$ is an isomorphism for every
$N\in{}^{\coring{C}}\mathcal{M}_{f}.$ Finally, let
$M\in\mathcal{M}^{\coring{C}}$ such that $M_{A}$ is flat. We have
\[
M\square_{\coring{C}}-\simeq
\hom{\coring{C}}{-}{M}_{f}(-)^{\ast}:{}^{\coring{C}}\mathcal{M}_{f}\rightarrow
\mathcal{M}_{k}.
\]
Then, $M_{\coring{C}}$ is coflat (resp. faithfully coflat) iff
$M\square_{\coring{C}}-:{}^{\coring{C}}\mathcal{M}_{f}\rightarrow
\mathcal{M}_{k}$ is exact (resp. faithfully exact) (by Lemma
\ref{12}) iff
$\hom{\coring{C}}{-}{M}_{f}:\mathcal{M}_{f}^{\coring{C}%
}\rightarrow\mathcal{M}_{k}$ is exact (resp. faithfully exact) iff
$M_{\coring{C}}$ is injective (resp. an injective cogenerator) (by
Lemma \ref{12}).
\end{proof}

\begin{corollary}\label{injectiveinjector}
Let $N\in{}^{\coring{C}}\mathcal{M}^{\coring{D}}$ be a bicomodule.
Suppose that $A$ is a QF ring and $\coring{D}$ has a duality. If
$N$ is injective
(resp. injective cogenerator) in $\mathcal{M}%
^{\coring{D}}$ such that $N_{B}$ is flat, then $N$ is an injector
(resp. an injector-cogenerator) as a right $\coring{D}$-comodule.
\end{corollary}

\begin{proof}
Let $X_{A}$ be an injective (resp. an injective cogenerator)
module. Since $A$ is a QF ring, $X_{A}$ is projective. We have
then the natural isomorphism
\[
(X\otimes_{A}N)\square_{\coring{D}}-\simeq
X\otimes_{A}(N\square_{\coring{D}}-)
:{}^{\coring{D}}\mathcal{M}\rightarrow \mathcal{M}_{k}.
\]
By Proposition \ref{13}, $N_{\coring{D}}$ and $X_{A}$ are coflat
(resp. faithfully coflat), and then $X\otimes_{A}N$ is coflat
(resp. faithfully coflat). Now, since $X\otimes_{A}N$ is a flat
right $B$-module, and by Proposition \ref{13} $X\otimes_{A}N$ is
injective (resp. injective cogenerator) in
$\mathcal{M}^{\coring{D}}$.
\end{proof}

 The last two results allow to improve our general statements in
 Section \ref{Frobeniusgeneral} for corings having a duality.

\begin{proposition}\label{adjointpairduality}
Suppose that $\coring{C}$ and $\coring{D}$ have a duality.
Consider the following statements
\begin{enumerate} \item $(-\cotensor{\coring{C}} X, -
\cotensor{\coring{D}}\Lambda)$ is an adjoint pair of functors,
with $_{A}X$ and $\Lambda_{A}$ flat;\item $\Lambda$ is
quasi-finite injective as a right ${\coring{C}}$-comodule, with
$_{A}X$ and $\Lambda_{A}$ flat and $X\simeq
\cohom{\coring{C}}{\Lambda}{\coring{C}}$ in
$^{\coring{C}}\mathcal{M}^{\coring{D}}$.
\end{enumerate}
We have (1) implies (2), and the converse is true if in particular
$B$ is a QF ring.
\end{proposition}

\begin{proof}
$(1)\Rightarrow(2)$ From Proposition \ref{13}, $\coring{D}$ is
injective in $\mathcal{M}^{\coring{D}}$. Since the functor
$-\cotensor{\coring{C}} X$ is exact, $\Lambda\simeq
\coring{D}\cotensor{\coring{D}} \Lambda$ is injective in
$\mathcal{M}^{\coring{C}}$.
 \\ $(2)\Rightarrow(1)$ Assume
that $B$ is a QF ring. From Corollary \ref{injectiveinjector},
$\Lambda$ is quasi-finite injector as a right
${\coring{C}}$-comodule, and Proposition \ref{19} achieves the
proof.
\end{proof}

\begin{theorem}\label{leftrightFrobeniusduality}
Suppose that $\coring{C}$ and $\coring{D}$ have a duality. Let
$X\in{}^{\coring{C}}\mathcal{M}%
^{\coring{D}}$ and
$\Lambda\in{}^{\coring{D}}\mathcal{M}^{\coring{C}}.$ The following
statements are equivalent \begin{enumerate} \item
$(-\cotensor{\coring{C}} X, - \cotensor{\coring{D}}\Lambda)$ is a
Frobenius pair, with $X_{B}$ and $\Lambda_{A}$ flat; \item
$(\Lambda \cotensor{\coring{C}} -, X \cotensor{\coring{D}}-)$ is a
Frobenius pair, with $_{A}X$ and $_{B}\Lambda$ flat; \item $X$ and
$\Lambda$ are quasi-finite injector on both sides, and $X\simeq
\cohom{\coring{C}}{\Lambda}{\coring{C}}$ in
$^{\coring{C}}\mathcal{M}^{\coring{D}}$ and $\Lambda\simeq
\cohom{\coring{D}}{X}{\coring{D}}$ in
$^{\coring{D}}\mathcal{M}^{\coring{C}}$.

In particular, if $A$ and $B$ are QF rings, then the above
statements are equivalent to \item $X$ and $\Lambda$ are
quasi-finite injective on both sides, and $X\simeq
\cohom{\coring{C}}{\Lambda}{\coring{C}}$ in
$^{\coring{C}}\mathcal{M}^{\coring{D}}$ and $\Lambda\simeq
\cohom{\coring{D}}{X}{\coring{D}}$ in
$^{\coring{D}}\mathcal{M}^{\coring{C}}$.
\end{enumerate}
Finally, suppose that $\coring{C}$ and $\coring{D}$ are
cosemisimple corings. Let
$X\in{}^{\coring{C}}\mathcal{M}%
^{\coring{D}}$ and
$\Lambda\in{}^{\coring{D}}\mathcal{M}^{\coring{C}}.$ The following
statements are equivalent \begin{enumerate} \item
$(-\cotensor{\coring{C}} X, - \cotensor{\coring{D}}\Lambda)$ is a
Frobenius pair; \item $(\Lambda \cotensor{\coring{C}} -, X
\cotensor{\coring{D}}-)$ is a Frobenius pair; \item $X$ and
$\Lambda$ are quasi-finite on both sides, and $X\simeq
\cohom{\coring{C}}{\Lambda}{\coring{C}}$ in
$^{\coring{C}}\mathcal{M}^{\coring{D}}$ and $\Lambda\simeq
\cohom{\coring{D}}{X}{\coring{D}}$ in
$^{\coring{D}}\mathcal{M}^{\coring{C}}$.
\end{enumerate}
\end{theorem}

\begin {proof}
At first we will prove the first part. In view of Theorem
\ref{leftrightFrobenius} and Theorem \ref{20a} it suffices to show
that if $(-\cotensor{\coring{C}} X, -
\cotensor{\coring{D}}\Lambda)$ is a Frobenius pair, the condition
``$X_{\coring{D}}$ and $\Lambda_{\coring{C}}$ are coflat'' is
equivalent to ``$X_{B}$ and $\Lambda_{A}$  are flat''. Indeed, the
first implication is obvious, for the converse, assume that
$X_{B}$ and $\Lambda_{A}$ are flat. By Proposition
\ref{adjointpairduality}, $X$ and $\Lambda$ are injective in
$\mathcal{M}^{\coring{D}}$ and $\mathcal{M}^{\coring{C}}$
respectively, and they are coflat by Proposition \ref{13}. The
particular case is straightforward from Proposition
\ref{adjointpairduality} and the above equivalences .
\\ Now we will show the second part. Since every comodule category
over a cosemisimple coring is a spectral category (see \cite[p.
128]{Stenstrom:1975}), and by Proposition \ref{13}, every comodule
(resp. every bicomodule) over a cosemisimple coring is coflat
(resp. injector) (we can see it directly by using the fact that
every additive functor between abelian categories preserves split
exactness). The use of Theorem \ref{leftrightFrobenius} finishes
then proof.
\end{proof}

\begin{remark}
\begin{enumerate}
\item The equivalence ``$(1)\Leftrightarrow(4)$'' of the last
theorem is a generalization of \cite[Theorem
3.3]{Castano/Gomez/Nastasescu:1999}. The proof of \cite[Theorem
3.3]{Castano/Gomez/Nastasescu:1999} give an alternative proof of
``$(1)\Leftrightarrow(4)$'' of Theorem
\ref{leftrightFrobeniusduality}, using Proposition \ref{11}. \item
The adjunction of Proposition \ref{19a} and Proposition
\ref{adjointpairduality} generalize the coalgebra version of
Morita's theorem \cite[Theorem
4.2]{Caenepeel/DeGroot/Militaru:2002}.
\end{enumerate}
\end{remark}

\begin{example}
Let $A$ be a $k$-algebra. Put $\coring{C}=A$ and $\coring{D}=k.$
The bicomodule $A\in{}^{\coring{C}}\mathcal{M}^{\coring{D}}$ is
quasi-finite as a right $\coring{D}$-comodule. $A$ is an injector
as a right $\coring{D} $-comodule if and only if the $k$-module
$A$ is flat. If we take $A=k=\mathbb{Z} $, the bicomodule $A$ is
quasi-finite and injector as a right $\coring{D}$-comodule but it
is not injective in $\mathcal{M}^{\coring{D}}.$ Hence, the
assertion ``$(-\cotensor{\coring{C}} X, -
\cotensor{\coring{D}}\Lambda)$ is an adjoint pair of functors''
does not implies in general the assertion ``$\Lambda$ is
quasi-finite injective as a right ${\coring{C}}$-comodule and
$X\simeq \cohom{\coring{C}}{\Lambda}{\coring{C}}$ in
$^{\coring{C}}\mathcal{M}^{\coring{D}}$'', and the following
statements are not equivalent in general:
\begin{enumerate} \item $(-\cotensor{\coring{C}} X, -
\cotensor{\coring{D}}\Lambda)$ is a Frobenius pair; \item $X$ and
$\Lambda$ are quasi-finite injective on both sides, and $X\simeq
\cohom{\coring{C}}{\Lambda}{\coring{C}}$ in
$^{\coring{C}}\mathcal{M}^{\coring{D}}$ and $\Lambda\simeq
\cohom{\coring{D}}{X}{\coring{D}}$ in
$^{\coring{D}}\mathcal{M}^{\coring{C}}$.
\end{enumerate}
On the other hand, there exists a commutative self-injective ring
wich is not coherent. By a theorem of S.U.~Chase (see for example
\cite[Theorem 19.20]{Anderson/Fuller:1992}), there exists then a
$k$-algebra $A$ which is injective, but not flat as $k$-module.
Hence, the bicomodule
$A\in{}^{\coring{C}}\mathcal{M}^{\coring{D}}$ is quasi-finite and
injective as a right $\coring{D}$-comodule, but not an injector as
a right $\coring{D}$-comodule.
\end{example}

\section{Applications to induction functors}\label{Frobmor}

We start this section by recalling from \cite{Gomez:2002}, that a
coring homomorphism from the coring $\coring{C}$\ to the coring
$\coring{D}$ is a pair $(\varphi,\rho)$, where $\rho:A\rightarrow
B$ is a homomorphism of $k$-algebras and
$\varphi:\coring{C}\rightarrow\coring{D}$ is a homomorphism of
$A$-bimodules such that
\[
\epsilon_{\coring{D}}\circ\varphi=\rho\circ\epsilon_{\coring{C}}\qquad\textrm{and
}\qquad\Delta_{\coring{D}}\circ\varphi=\omega_{\coring{D},\coring{D}}\circ(
\varphi\otimes_{A}\varphi)\circ\Delta_{\coring{C}},
\]
where $\omega_{\coring{D},\coring{D}}:\coring{D}\otimes_{A}\coring{D}%
\rightarrow\coring{D}\otimes_{B}\coring{D}$ is the canonical map
induced by $\rho:A\rightarrow B.$

\medskip

Now we will characterize when the induction functor $ - \tensor{A}
B : \rcomod{\coring{C}} \rightarrow \rcomod{\coring{D}}$ defined
in \cite[Proposition 5.3]{Gomez:2002} is a Frobenius functor. The
coaction of $\coring{D}$ over $M \tensor{A} B$ is given, when
expressed in Sweedler's sigma notation, by
\[
\rho_{M \tensor{A} B}(m \tensor{A}b) = \sum m_{(0)} \tensor{A} 1
\tensor{B} \varphi(m_{(1)})b,
\]
where $M$ is a right $\coring{C}$--comodule with coaction
$\rho_M(m) = \sum m_{(0)} \tensor{A} m_{(1)}$.

\begin{theorem}\label{23}
Let $(\varphi,\rho) :\coring{C}\rightarrow\coring{D}$ be a
homomorphism of corings such that $_{A}\coring{C}$ and
$_{B}\coring{D} $ are flat. The following statements are
equivalent\begin{enumerate}[(a)] \item
$-\otimes_{A}B:\mathcal{M}^{\coring{C}}\rightarrow\mathcal{M}%
^{\coring{D}}$ is a Frobenius functor;\item the $\coring{C}%
-\coring{D}$-bicomodule $\coring{C}\otimes_{A}B$ is quasi-finite
and injector as a right $\coring{D}$-comodule and there exists an
isomorphism of $\coring{D}-\coring{C}$-bicomodules
$\cohom{\coring{D}}{\coring{C}\otimes_{A}B}{\coring{D}}\simeq
B\otimes _{A}\coring{C}.$
\end{enumerate}
Moreover, if $\coring{C}$ and $\coring{D}$ are coseparable, then
the condition ``injector'' in $(b)$ can be deleted.
\end{theorem}

\begin{proof}
First observe that $-\otimes_{A}B$ is a Frobenius functor if and
only if $(-\otimes_{A}B,-\square_{\coring{D}}(B\otimes
_{A}\coring{C}))$ is a Frobenius pair (by \cite[Proposition
5.4]{Gomez:2002}). A straightforward computation shows that the
map $\rho_M \tensor{A} B : M \tensor{A} B \rightarrow M \tensor{A}
\coring{C} \tensor{A} B$ is a homomorphism of
$\coring{D}$--comodules. We have thus a commutative diagram in
$\rcomod{\coring{C}}$ with exact row
\[
\xymatrix{ 0 \ar[r] & M \cotensor{\coring{C}} (\coring{C}
\tensor{A} B) \ar^{\iota}[r] & M \tensor{A} \coring{C} \tensor{A}
B \ar^-{\omega_{M,\coring{C} \tensor{A} B}}[rr] & & M \tensor{A}
\coring{C} \tensor{A} \coring{C} \tensor{A} B \\
 & M \tensor{A} B \ar^{\psi_M}[u] \ar_-{\rho_{M \tensor{A} B}}[ur] & &
 &},
\]
where $\psi_M$ is defined by the universal property of the kernel.
Since ${}_B\coring{D}$ is flat, to prove that $\psi_M$ is an
isomorphism of $\coring{D}$--comodules it is enough to check that
it is bijective, as the forgetful functor $U : \rcomod{\coring{D}}
\rightarrow \rmod{B}$ is faithfully exact. Some easy computations
show that the map $(M \tensor{A} \epsilon_{\coring{C}} \tensor{A}
B)\circ \iota$ is the inverse in $\rmod{B}$ to $\psi_M$. From
this, we deduce a natural isomorphism $ \psi : -\otimes_{A}B\simeq
-\square_{\coring{C}}( \coring{C}\otimes_{A}B)$. The equivalence
between (a) and (b) is then obvious from Proposition \ref{19} and
Proposition \ref{19a}.
\end{proof}

When applied to the case where $\coring{C} = A$ and $\coring{D} =
B$ are the trivial corings, Theorem \ref{23} gives functorial
Morita's characterization of Frobenius ring extensions
\cite{Morita:1965}. It is then reasonable to give the following
definition.

\begin{definition}
Let $(\varphi,\rho) :\coring{C}\rightarrow\coring{D}$ be a
homomorphism of corings such that $_{A}\coring{C}$ and
$_{B}\coring{D} $ are flat. It is said to be a \emph{right
Frobenius} morphism of corings if $ - \tensor{A} B :
\rcomod{\coring{C}} \rightarrow \rcomod{\coring{D}}$ is a
Frobenius functor.
\end{definition}

 The following generalize \cite[Theorem
3.5]{Castano/Gomez/Nastasescu:1999}.

\begin{theorem}\label{24}
Suppose that the algebras $A$ and $B$ are QF rings.\newline Let
$(\varphi,\rho):\coring{C}\rightarrow\coring{D}$ be a homomorphism
of corings such that the modules $_{A}\coring{C}$,
$_{B}\coring{D}$ and $\coring{D}_{B}$ are projective. Then the
following statements are equivalent
\begin{enumerate}[(a)]
\item $-\otimes_{A}B:\mathcal{M}^{\coring{C}}\rightarrow
\mathcal{M}^{\coring{D}}$ is a Frobenius functor;\item the
$\coring{C}-\coring{D}$-bicomodule $\coring{C}\otimes_{A}B$ is
quasi-finite as a right $\coring{D}$-comodule,
$(\coring{C}\otimes_{A}B) _{\coring{D}}$ is injective and there
exists an isomorphism of
$\coring{D}%
-\coring{C}$-bicomodules
$\cohom{\coring{D}}{\coring{C}\otimes_{A}B}{\coring{D}}\simeq
B\otimes_{A}\coring{C}.$
\end{enumerate}
\end{theorem}

\begin{proof}
Obvious from Proposition \ref{adjointpairduality}.
\end{proof}

Now, suppose that the forgetful functor $\mathcal{M}^{\coring{C}%
}\rightarrow\mathcal{M}_{A}$ is a Frobenius functor. Then the
functor
$-\otimes_{A}\coring{C}:\mathcal{M}_{A}\rightarrow\mathcal{M}_{A}$
is also a Frobenius functor (since it is a composition of two
Frobenius functors) and $_{A}\coring{C}$ is finitely generated
projective. On the other hand, since
$-\otimes_{A}\coring{C}:\mathcal{M}_{A}\rightarrow\mathcal{M}%
^{\coring{C}}$ is a left adjoint to
$\hom{\coring{C}}{\coring{C}}{-}
:\mathcal{M}^{\coring{C}}\rightarrow\mathcal{M}_{A}.$ Then
$\hom{\coring{C}}{\coring{C}}{-}$ is a Frobenius functor.
Therefore, $\coring{C}$ is finitely generated projective in
$\mathcal{M}^{\coring{C}},$ and hence in $\mathcal{M}_{A}.$

\begin{lemma}\label{25}
Let $R$ be the opposite algebra of $^{\ast}\coring{C}$.
\begin{enumerate}[(1)]
\item $\coring{C}\in{}^{\coring{C}}\mathcal{M}^{A}$ is
quasi-finite (resp. quasi-finite and injector) as a right
$A$-comodule if and only if $_{A}\coring{C}$ is finitely generated
projective (resp. $_{A}\coring{C}$ is finitely generated
projective and $_{A}R$ is flat). Let $\cohom{A}{\coring{C}}{-}
=-\otimes_{A}R:\mathcal{M}^{A}\rightarrow\mathcal{M}^{\coring{C}}$
be the co-hom functor.\item If $_{A}\coring{C}$ is finitely
generated projective and $_{A}R$ is flat, then
\[
_{A}\cohom{A}{\coring{C}}{A}_{\coring{C}}\simeq{}_{A}%
R_{\coring{C}},
\]
where the right $\coring{C}$-comodule structure of $R$ is defined
as in \cite[Lemma 4.3]{Brzezinski:2002}.
\end{enumerate}
\end{lemma}

\begin{proof}
(1) Straightforward from \cite[Example 4.3]{Gomez:2002}.\newline
(2) From \cite[Lemma 4.3]{Brzezinski:2002}, the forgetful functor
$\mathcal{M}^{\coring{C}}\rightarrow\mathcal{M}%
_{A}$ is the composition of functors
$\mathcal{M}^{\coring{C}}\rightarrow
\mathcal{M}_{R}\rightarrow\mathcal{M}_{A}.$ By \cite[Proposition
4.2]{Gomez:2002}, $\cohom{A}{\coring{C}}{-}$ is a left adjoint to
$-\square_{\coring{C}}\coring{C}:\mathcal{M}^{\coring{C}}\rightarrow
\mathcal{M}_{A}$ which is isomorphic to the forgetful functor
$\mathcal{M}%
^{\coring{C}}\rightarrow\mathcal{M}_{A}.$ Then
$\cohom{A}{\coring{C}}{-}$ is isomorphic to the composition of
functors
\[
\xymatrix@1{\mathcal{M}_{A}\ar[r]^{-\otimes
_{A}R}&\mathcal{M}_{R}\ar[r]& \mathcal{M}^{\coring{C}}.}
\]
In particular, $_{A}\cohom{A}{\coring{C}}{A}_{\coring{C}}\simeq
{}_{A}( A\otimes_{A}R)_{\coring{C}}\simeq{}_{A}R_{\coring{C}}.$
\end{proof}

\begin{corollary}\label{26}
(\cite[27.10]{Brzezinski/Wisbauer:2003}) \newline Let $\coring{C}$
be an $A$-coring and let $R$ be the opposite algebra of
$^{\ast}\coring{C}.$ Then the following statements are equivalent
\begin{enumerate}[(a)]
\item The forgetful functor $F:\mathcal{M}^{\coring{C}%
}\rightarrow\mathcal{M}_{A}$ is a Frobenius functor;\item
$_{A}\coring{C}$ is finitely generated projective and
$\coring{C}\simeq R$ as $(A,R)$-bimodules, where $\coring{C}$ is a
right $R$-module via $c.r=c_{(1)}.r(c_{(2)})$, for all
$c\in\coring{C}$ and $r\in R.$
\end{enumerate}
\end{corollary}

\begin{proof}
Straightforward from Theorem \ref{23} and Lemma \ref{25}.
\end{proof}

The following proposition gives sufficient conditions to have that
a morphism of corings is right Frobenius if and only if it is left
Frobenius. Note that it says in particular that the notion of
Frobenius homomorphism of coalgebras over fields (by (b)) or of
rings (by (d)) is independent on the side. Of course, the latter
is well-known.

\begin{proposition}\label{27}
Let $(\varphi,\rho):\coring{C}\rightarrow\coring{D}$ be a
homomorphism of corings such that $_{A}\coring{C}$,
$_{B}\coring{D}$, $\coring{C}_{A}$ and $\coring{D}_{B}$ are flat.
Assume that at least one of the following holds
\begin{enumerate}[(a)] \item $\coring{C}$ and $\coring{D}$ have a
duality, and $_{A}B$ and $B_{A}$ are flat;\item $A$ and $B$ are
von Neumann regular
rings;\item $B\otimes_{A}\coring{C}$ is coflat in $^{\coring{D}%
}\mathcal{M}$ and $\coring{C}\otimes_{A}B$ is coflat in $\mathcal{M}%
^{\coring{D}}$ and $_{A}B $ and $B_{A}$ are flat;\item
$\coring{C}$ and $\coring{D}$ are coseparable corings.
\end{enumerate} Then the following statements are equivalent
\begin{enumerate}
\item $-\otimes_{A}B:\mathcal{M}^{\coring{C}%
}\rightarrow\mathcal{M}^{\coring{D}}$ is a Frobenius functor;\item
$B\otimes_{A}-:{}^{\coring{C}}\mathcal{M}\rightarrow{}^{\coring{D}%
}\mathcal{M}$ is a Frobenius functor.
\end{enumerate}
\end{proposition}

\begin{proof}
Obvious from Theorem \ref{leftrightFrobenius} and Theorem
\ref{leftrightFrobeniusduality}.
\end{proof}

Let us finally show to derive from our results a remarkable
characterization of the so called \emph{Frobenius corings}.

\begin{corollary}\label{28}
(\cite[27.8]{Brzezinski/Wisbauer:2003})\newline The following
statements are equivalent \begin{enumerate}[(a)] \item the
forgetful functor
$\mathcal{M}^{\coring{C}}\rightarrow\mathcal{M}_{A}$ is
a Frobenius functor;\item the forgetful functor $^{\coring{C}%
}\mathcal{M}\rightarrow{}_{A}\mathcal{M}$ is a Frobenius
functor;\item there exist an $(A,A) $-bimodule map
$\eta:A\rightarrow \coring{C}$ and a $(\coring{C},\coring{C})
$-bicomodule map
$\pi:\coring{C}\otimes_{A}\coring{C}\rightarrow\coring{C}$ such
that $\pi(\coring{C}\otimes_{A}\eta) =\coring{C}=\pi(
\eta\otimes_{A}\coring{C}) .$
\end{enumerate}
\end{corollary}

\begin{proof}
The proof of ``$(1)\Leftrightarrow(4)$'' in Proposition \ref{19a}
for $X=\coring{C}\in{}^{A}\mathcal{M}^{\coring{C}}$ and
$\Lambda=\coring{C}\in{}^{\coring{C}}\mathcal{M}^{A}$ remains
valid for our situation. It remains to see that the condition (4)
in this case is exactly the condition (c).
\end{proof}

\section{Applications to entwined modules}

In this section we particularize some our results in Section
\ref{Frobmor} to the category of entwined modules. We adopt the
notations of \cite{Caenepeel/Militaru/Zhu:2002}. We start with
some remarks.

\begin{enumerate}[(1)]

\item  Let a right-right entwining structure
$(A,C,\psi)\in\mathbb{E}_{\bullet}^{\bullet}(k)$
 and a left-left entwining structure $(B,D,\varphi)\in{}_{\bullet
 }^{\bullet}\mathbb{E}(k)$. The category of two-sided entwined modules
$_{B}^{D}\mathcal{M}(\varphi,\psi)_{A}^{C}$ defined in \cite[pp.
68--69]{Caenepeel/Militaru/Zhu:2002} is isomorphic to the category
of bicomodules $^{D\otimes B}\mathcal{M}^{A\otimes C}$ over the
associated corings.

\item If $(A,C,\psi)$ and $(A^{\prime },C^{\prime},\psi^{\prime})$
belong to $\mathbb{E}_{\bullet}^{\bullet }(k)$ such that $\psi$ is
an isomorphism, then $\psi$ is an isomorphism of corings (see
\cite[Proposition 34]{Caenepeel/Militaru/Zhu:2002}), and
consequently if the coalgebra $C$ is flat as a $k$-module, then
the modules $_{A}(A\otimes C)$ and $(A\otimes C)_{A}$ are flat,
and \[
^{A\otimes C}\mathcal{M}^{A^{\prime}\otimes C^{\prime}}\simeq{}_{A}%
^{C}\mathcal{M}(\psi^{-1},\psi^{\prime})_{A^{\prime}}^{C^{\prime}}.
\]

\item  Let $(\alpha,\gamma):(A,C,\psi)\rightarrow(A^{\prime
},C^{\prime},\psi^{\prime})$ be a morphism in
$\mathbb{E}_{\bullet}^{\bullet }(k)$. We know that
$(\alpha\otimes\gamma,\alpha):A\otimes{C}\rightarrow{A^{\prime}}\otimes{C^{\prime}}$
is a morphism of corings. The functor $F$ defined in \cite[Lemma
8]{Caenepeel/Militaru/Zhu:2002} satisfies the commutativity of the
diagram
\[
\xymatrix{\mathcal{M}^{A\otimes
C}\ar[rr]^{-\otimes_{A}A^{\prime}}\ar[d]_\simeq & &
\mathcal{M}^{A^{\prime}\otimes C^{\prime}}\ar[d]^\simeq \\
\mathcal{M}(\psi)_{A}^{C}\ar[rr]_{-\otimes_{A}A^{\prime}} & &
\mathcal{M}(\psi^{\prime})_{A^{\prime}}^{C^{\prime}},}
\]
where $-\otimes_{A}A^{\prime}:\mathcal{M}^{A\otimes C}\rightarrow
\mathcal{M}^{A^{\prime}\otimes C^{\prime}}$ is the induction
functor defined in \cite[Proposition 5.3]{Gomez:2002}.

\end{enumerate}

We have the following result concerning the category of entwined
modules:

\begin{theorem}\label{entwined modules}
Let $(\alpha,\gamma):(A,C,\psi)\rightarrow(A^{\prime
},C^{\prime},\psi^{\prime})$ be a morphism in
$\mathbb{E}_{\bullet}^{\bullet }(k)$, such that $_{k}C$ and
$_{k}D$ are flat.
\begin{enumerate} \item The following statements are equivalent\begin{enumerate}[(a)]
\item The functor
$F=-\otimes_{A}A^{\prime}:\mathcal{M}(\psi)_{A}^{C}\rightarrow
\mathcal{M}(\psi^{\prime})_{A^{\prime}}^{C^{\prime}}$ defined in
\cite[Lemma 8]{Caenepeel/Militaru/Zhu:2002} is a Frobenius
functor; \item the $ A\otimes{C}-A^{\prime}\otimes
C^{\prime}$-bicomodule $(A\otimes C)\otimes_{A}A^{\prime}$ is
quasi-finite injector as a right $A^{\prime}\otimes
C^{\prime}$-comodule and there exists an isomorphism of
$A^{\prime}\otimes C^{\prime}-A\otimes C$-bicomodules
$\cohom{A^{\prime}\otimes
C^{\prime}}{(A\otimes{C})\otimes_{A}A^{\prime}}{A^{\prime}\otimes
C^{\prime}}\simeq A^{\prime}\otimes_{A}(A\otimes{C}).$
\end{enumerate}
Moreover, if $A\otimes C$ and $A^{\prime}\otimes C^{\prime}$ are
coseparable corings, then the condition ``injector'' in (b) can be
deleted.

 \item If $A$ and $B$ are QF rings and the module
$(A^{\prime}\otimes C^{\prime})_{A^{\prime}}$ is projective, then
the following are equivalent\begin{enumerate}[(a)] \item The
functor
$F=-\otimes_{A}A^{\prime}:\mathcal{M}(\psi)_{A}^{C}\rightarrow\mathcal{M}
(\psi^{\prime})_{A^{\prime}}^{C^{\prime}}$ defined in \cite[Lemma
8]{Caenepeel/Militaru/Zhu:2002} is a Frobenius functor; \item the
$ A\otimes C-A^{\prime}\otimes C^{\prime}$-bicomodule
$(A\otimes{C})\otimes_{A}A^{\prime}$ is quasi-finite and injective
as a right $A^{\prime}\otimes C^{\prime}$-comodule and there
exists an isomorphism of $A^{\prime}\otimes C^{\prime}-A\otimes
C$-bicomodules $\cohom{A^{\prime}\otimes
C^{\prime}}{(A\otimes{C})\otimes_{A}A^{\prime}}{A^{\prime}\otimes
C^{\prime}}\simeq A^{\prime}\otimes_{A}(A\otimes{C}).$
\end{enumerate}
\end{enumerate}
\end{theorem}

\begin{proof}
Follows from Theorem \ref{23} and Theorem \ref{24}.
\end{proof}

\begin{remark}
Let a right-right entwining structure
$(A,C,\psi)\in\mathbb{E}_{\bullet}^{\bullet}(k)$. The
coseparability of the coring ${A}\tensor{}{C}$ is characterized in
\cite[Theorem 38(1)]{Caenepeel/Militaru/Zhu:2002} (see also
\cite[Corollary 3.6]{Brzezinski:2002}).
\end{remark}

\section*{Acknowledgements} We would like to thank the Professor
Edgar Enochs to have communicated to us the example of a
commutative self-injective ring which is not coherent.

\end{document}